\documentclass[reqno,12pt]{amsart}
\usepackage{amsmath}
\usepackage{amssymb}
\usepackage{amsthm}
\newtheorem{theorem}{Theorem} % [section]
\newtheorem{remark}{Remark}
\newtheorem{lemma}{Lemma}

\begin{document}

\title [Estimates for the rate of strong approximation]{Estimates for the rate of strong approximation in
Hilbert space}

\author[F. ~G\"otze]{Friedrich G\"otze$^1$}
\author [A.Yu. Zaitsev]{Andrei Yu. Zaitsev$^{1,2}$}

\address{Fakult\"at f\"ur Mathematik,\smallskip
Universit\"at Bielefeld, Postfach 100131, D-33501 Bielefeld,
Germany}
\email{goetze@math.uni-bielefeld.de}

\address{St.~Petersburg Department of Steklov Mathematical Institute,
\\
Fontanka 27, St.~Petersburg 191023, Russia}
\email{zaitsev\@pdmi.ras.ru}



\parindent.6cm

\date{FMarch 18, 2011}

\maketitle

\footnotetext[1]{Research supported by the SFB 701 in Bielefeld and by grant
RFBR-DFG 09-01-91331.}
%\newline\indent
\footnotetext[2]{Research supported by grant 10-01-00242,
by grant of leading Scientific scools 4472.2010.1 and by a program of
fundamental researches of Russian Academy of Sciences ``Modern
problems of fundamental mathematics''.
}

\section{Introduction}

The aim of this paper is to investigate, which infinite
dimensional consequences follow from the main results of recently
published paper of the authors~\cite{GZ2} (see Theorems~\ref{th3}
and \ref{th4}). We show that the finite dimensional Theorem
\ref{th4} implies meaningful estimates for the rate of strong
Gaussian approximation of sums of i.i.d. Hilbert space valued
random vectors $\xi_j$
 with finite
moments ${\mathbf E}\,\left\| \xi _j\right\|^\gamma$, $\gamma>2$.
We show that the rate of approximation depends substantially on
the rate of decay of the sequence of eigenvalues of the covariance
operator of summands.

Below we need some notation.  The distribution of a random vector
$\xi $ will be denoted by $\mathcal L(\xi )$. The corresponding
covariance operator will be denoted by $\hbox{\rm cov}\,\xi$. We
denote $\log ^\ast b=\max \left\{ 1,\log b\right\} $ for $b>0$.
 We shall write $A\ll_t B$, if there exists a
positive quantity $c(t)$ depending only on $t$ and such that
$A\leqslant c(t)\, B$. We shall also write $A\asymp_t B$, if
$A\ll_t B\ll_t A$. The absence of lower indices means that the
corresponding constants are absolute.

We consider the following well-known problem. Let $\xi _1,\dots
,\xi _n$ be independent random vectors with zero means and finite
moments of second order.  One has to construct on the same
probability space a sequence of independent random vectors
$X_1,\dots ,X_n$ and  independent Gaussian random vectors
$Y_1,\dots ,Y_n$ such that
$$
\mathcal L(X_j)=\mathcal L(\xi_j), \quad{\mathbf E}\,Y_j=0,\quad
\hbox{\rm cov}\,Y_j=\hbox{\rm cov}\,X_j, \quad\quad j=1,\dots ,n,
$$
 and the quantity
\begin{equation}
\Delta_n (X,Y)=\max_{1\leqslant s\leqslant n}\;\Bigl\|\,\sum_{j=1}^sX_j-
\sum_{j=1}^sY_j\,\Bigr\| \label{1.1}
\end{equation} would be as small as possible with sufficiently large
probability.
The estimation of the rate of strong approximation in
the invariance principle may be reduced just to this problem. We omit
the detailed history of the problem referring the reader to
G\"otze and Zaitsev~\cite{5} and Zaitsev~\cite{13}.

For brevity, instead of writting out the properties of the vectors
$X_1,\dots ,X_n$  and $Y_1,\dots ,Y_n$ listed above we simply say
that {\it there exists a construction} having additional
properties which are mentioned explicitly in the text. As a rule,
we consider the case where the vectors $\xi _1,\dots ,\xi _n$ are
identically distributed with some random vector~$Z$ and, in
conditions of Theorems, we mention just this vector.

In this paper,
 we obtain infinite dimensional analogues of the following result of
Sakhanenko \cite{12} in the case of i.i.d. summands.

\begin{theorem}\label{th2} Let $\xi _1,\xi_2,\dots ,\xi
_n$ be independent random variable with ${{\mathbf E}\,\xi _j=0}$,
$j=1,\dots ,n$. Let $\gamma > 2$ and
$$
L_\gamma =\sum_{j=1}^n{\mathbf E}\,\left\vert \xi _j\right\vert
^\gamma <\infty.
$$
Then there exists a construction such that
\begin{equation} {\mathbf E}\,\bigl( \Delta_n (X,Y)\bigr)^\gamma \ll
_\gamma L_\gamma .\label{1.2} \end{equation} \end{theorem}

\medskip

It should be mentioned that, in Sakhanenko \cite{12}, one can find more
general results. In Sakhanenko \cite{12}, it is observed that
inequality~\eqref{1.2} implies the well-known Rosenthal
inequality (\cite{10},
\cite{11}, see Lemma~\ref{lr1}).

Upon the natural normalization, we see that \eqref{1.2} is
equivalent to
$$
{\mathbf E}\,\bigl( \Delta_n (X,Y)/\sigma\bigr)^\gamma \ll_
\gamma L_\gamma /\sigma^\gamma,
$$
where $\sigma^2= \hbox{\rm Var}\bigl(\,\sum_{j=1}^n\xi _j\bigr)$.
 It is clear that $L_\gamma /\sigma^\gamma$, $2<\gamma\leqslant3$, is
the well-known Lyapunov fraction involved in the Lyapunov and
Ess\'een bounds for the Kolmogorov distance in the CLT.

In this paper, we prove Theorems \ref{th6} and \ref{th9} which are
quite elementary consequences of Theorem~\ref{th3}, proved by the
authors in \cite{5} and~\cite{GZ2}. In~\cite{5}, we consider the
case of independent and (in general) non-identically distributed
summands. Theorem~\ref{th3} shows what follows from the results of
\cite{5} in a particular case, where summands are identically
distributed. Theorem~\ref{th3} is a multidimensional version of
Theorem~\ref{th2} for identically distributed summands.
\medskip

Denote  by $ \mathbf H$ the separable Hilbert space, which
consists of all real sequences~${x=(x_1,x_2,\dots )}$, for which
${\|x\|^2= x_1^2+x_2^2+\dots <\infty}$. Also put $\|x\|_\infty=
\max_j|x_j|$, ${x^{(d)}=(x_1,x_2,\dots,x_d )}\in\mathbf R^d$ and
$${x^{[d]}=(0,0,\dots,0, x_d,x_{d+1},\dots )}\in\mathbf H.$$

The formulations of our results involve a random vector
 $$Z=(Z_1,Z_2,\dots),$$ taking values in $ \mathbf H$ or $ \mathbf R^d$.
Independent copies of the vector $Z$ are to be constructed on the
same probability space with a sequence of independent Gaussian
random vectors. Without loss of generality, we assume, that the
coordinates of the vector~$Z$ are uncorrelated, and
\begin{equation}\label{in443}
\sigma_1^2\geqslant\sigma_2^2\geqslant\dots\geqslant\sigma_m^2\geqslant\dots,\quad
\hbox{where} \quad\sigma_m^2=\mathbf E\,Z_m^2,\quad\quad
m=1,2,\dots,
\end{equation}
 and
\begin{equation}\label{in44}\mathbb D=\hbox{\rm cov}\,Z,\quad
\mathbb D_d=\hbox{\rm cov}\,Z^{(d)},\quad
B_d^2=\sum_{m=d+1}^\infty\sigma_m^2={\mathbf
E}\,\bigl\|Z^{[d]}\bigr\|^2.
\end{equation}
 In particular,
\begin{equation}\label{in45}B_0^2=\sum_{m=1}^\infty\sigma_m^2={\mathbf E}\,\bigl\|Z\bigr\|^2.
\end{equation}

Moreover, in the formulations of our results, a  parameter $\psi$
satisfying
\begin{equation}
21/2<\psi\leqslant11.\label{1.111}
\end{equation} is involved.
In the sequel, Many constants below depend on~$\psi$. In order to
avoid this complication, one can simply take $\psi=11$.

\begin{theorem}\label{th3} Let $\psi$ satisfy \eqref{1.111} and let \,$Z$
be an $ \mathbf R^d$-valued random vector with $\sigma_d^2>0$,
${\mathbf E}\,Z=0$ and\/ ${\mathbf E}\,\left\Vert Z\right\Vert
^\gamma <\infty$, for some $\gamma \geqslant 2$. Then there exists
a construction such that
\begin{equation}
{\mathbf E}\, \bigl(\Delta_n (X, Y)\bigr)^\gamma\ll_{\psi,\gamma}
A\,\bigl( \sigma_{1}/\sigma_{d}\bigr)^\gamma \, n\, {\mathbf E}\,
\|Z\|^\gamma,\quad\text{for all } n,\label{1.10}
\end{equation}
where
\begin{equation}
A=A(\gamma,\psi,d)=\max\Big\{d^{\psi\gamma},\;d^{\frac{\gamma(\gamma+2)}4}\,
(\log^*d)^{\frac{\gamma(\gamma+1)}2}\Big\}.\label{1.11s}
\end{equation}\end{theorem}\medskip

We need a slightly different version of the finite dimensional
result. The following statement is proved in \cite{GZ2} while
proving Theorem~\ref{th3}.

\begin{theorem} \label{th4} Let $\psi$ satisfy \eqref{1.111} and let \,$Z$
be an $ \mathbf R^d$-valued random vector with $\sigma_d^2>0$,
${\mathbf E}\,Z=0$ and\/ ${\mathbf E}\,\left\Vert  Z\right\Vert
^\gamma <\infty$, for some $\gamma \geqslant 2$. There exists a
positive quantity $c_1(\gamma)$ depending only on $\gamma$ and
such that if
\begin{equation}
 \label{1.10bt}
 C(\gamma)\, d^{\gamma/2}(\log^*
d)^{\gamma+1}\,\bigl(\mathbf E\left\|\mathbb
D^{-1/2}Z\right\|^\gamma\bigr)^{2/\gamma}\leqslant n^{1-2/\gamma},
\end{equation}for some positive integer~$n$,
then there exists a construction such that
\begin{equation}
{\mathbf E}\, \bigl(\Delta_n (\mathbb D^{-1/2}X, \mathbb
D^{-1/2}Y)\bigr)^\gamma\ll_{\gamma,\psi}d^{\psi\gamma}\, n\,
{\mathbf E}\, \|\mathbb D^{-1/2}Z\|^\gamma.\label{1.10b}
\end{equation} \end{theorem}\medskip

\begin{remark}\rm  In \cite{GZ2}, the statements of Theorems~\ref{th3} and~\ref{th4}
involve, for $d^{\psi\gamma}$, the additional logarithmic factor
$(\log ^\ast d)^{2\gamma}$. We can easily eliminate it, observing
that we allow the constants in \eqref{1.10} and \eqref{1.10b} to
depend on $\psi$ satisfying \eqref{1.111}.
\end{remark}

\begin{remark}\rm
If condition \eqref{1.10bt} is not satisfied, then the estimates
in Theorem~\ref{th3} are obtained not due to a successful
  approximation,
but by estimating
${\mathbf E}\,\max_{1\leqslant s\leqslant n}\;\Bigl\|\,\sum_{j=1}^sX_j\,\Bigr\|^\gamma$ and ${\mathbf E}\,\max_{1\leqslant s\leqslant n}\;\Bigl\|\,
\sum_{j=1}^sY_j\,\Bigr\|^\gamma$ from above with the help of Lemma \ref{lr1} and inequality \eqref{1.16}.
Thus, the presence of condition \eqref{1.10bt} in the formulation of Theorem~\ref{th4}
 do not lead to a loss of information on the closeness of distributions in comparison with Theorem~\ref{th3}.
\end{remark}

The main results of this paper are Theorems~\ref{th6} and
\ref{th9}.

\begin{theorem}\label{th6} Let $\psi$ satisfy \eqref{1.111} and let \;$Z$
be an $ \mathbf H$-valued random vector with ${\mathbf E}\,Z=0$
and ${\mathbf E}\,\left\Vert  Z\right\Vert ^\gamma <\infty$, for
some $\gamma > 2$. If, for some fixed positive integers~$d$
and~$n$, the inequality
\begin{equation}\label{1.10ru}
C(\gamma)\, d^{\gamma/2}(\log^* d)^{\gamma+1}\,\bigl(\mathbf
E\,\bigl\|\mathbb
D_d^{-1/2}Z^{(d)}\bigr\|^\gamma\bigr)^{2/\gamma}\leqslant
n^{1-2/\gamma}
\end{equation} is valid, where  $C(\gamma)$ is defined in
Theorem~$\ref{th4}$,
then there exists a construction such that
\begin{equation}\label{1.10r}
{\mathbf E}\, \bigl(\Delta_n (X, Y)\bigr)^\gamma\ll_{\gamma
,\psi}d^{\psi\gamma}\, n\, \sigma_{1}^\gamma \, {\mathbf E}\,
\bigl\|\mathbb D_d^{-1/2}Z^{(d)}\bigr\|^\gamma+ \;n \,{\mathbf
E}\, \bigl\|Z^{[d]}\bigr\|^\gamma+
 (n\,B_d^2)^{\gamma/2}.\end{equation}
\end{theorem}\medskip

\begin{theorem}\label{th9} Let $\psi$ satisfy \eqref{1.111} and let \;$Z$
be a $ \mathbf H$-valued random vector with ${\mathbf E}\,Z=0$ and
${\mathbf E}\,\left\Vert  Z\right\Vert ^\gamma <\infty$, for some
$\gamma > 2$. If, for some fixed positive integers~$d$ and~$n$,
the inequality
\begin{equation}\label{1.10hu}
 C(\gamma)\,
d^{\gamma/2}(\log^*  d)^{\gamma+1}\,\bigl(\mathbf
E\left\|Z\right\|^\gamma\bigr)^{2/\gamma}\leqslant
n^{1-2/\gamma}\,\sigma_d^2,
\end{equation} is valid, where  $C(\gamma)$ is defined in
Theorem~$\ref{th4}$, then there exists a construction such that
\begin{equation}\label{1.10hr}
{\mathbf E}\, \bigl(\Delta_n (X, Y)\bigr)^\gamma\ll_{\gamma ,\psi}
d^{\psi\gamma}\, \bigl( \sigma_{1}/\sigma_{d}\bigr)^\gamma \, n\,
{\mathbf E}\, \|Z\|^\gamma+
 (n\,B_d^2)^{\gamma/2} .
\end{equation}
\end{theorem}\medskip

Theorems \ref{th6} and \ref{th9} make it possible to obtain
meaningful infinite dimensional estimates by a suitable choice of
dimension~$d$, for which the summands in the right-hand side of
inequality \eqref{1.10r} have approximately the same order in $n$.
Theorem \ref{th9} is an elementary consequence of Theorem
\ref{th6} and the inequality
\begin{equation}\label{bbn}
{\mathbf E}\, \bigl\|\mathbb
D_d^{-1/2}Z^{(d)}\bigr\|^\gamma\leqslant \sigma_{d}^{-\gamma} \,
{\mathbf E}\, \bigl\|Z^{(d)}\bigr\|^\gamma\leqslant
\sigma_{d}^{-\gamma} \, {\mathbf E}\, \|Z\|^\gamma.
\end{equation}

In general, Theorem \ref{th6} is sharper than Theorem \ref{th9}.
Many distributions with a regular behavior of moments satisfy the
relation
\begin{equation}\label{bn9}K=\sup_{1\leqslant d<\infty}d^{-\gamma/2}\,{\mathbf E}\, \|\mathbb
D_d^{-1/2}Z^{(d)}\|^\gamma<\infty ,\end{equation} which may lead
to a substantial improvement of the order of estimates. For
instance, if the vector $Z$ has independent coordinates $Z_m$,
then, by Lemma~\ref{lr2} of Section 2,
\begin{equation}\label{bb2}{\mathbf E}\, \|\mathbb
D_d^{-1/2}Z^{(d)}\|^\gamma= {\mathbf E}\,
\biggl(\sum\limits_{m=1}^d \frac{Z_m
^2}{\sigma_m^2}\biggr)^{\gamma/2} \ll_\gamma d^{\gamma/2}+
\sum\limits_{m=1}^d\sigma_m^{-\gamma} \,{\mathbf E}\,
|Z_m|^\gamma.
\end{equation}
Hence, $K<\infty$, if the sequence of moments $\sigma_m^{-\gamma}
\,{\mathbf E}\, |Z_m|^\gamma$ is bounded or grows not faster than
 $O(m^{(\gamma-2)/2})$. Observe that Lyapunov's inequality yields
 ${\mathbf E}\, \|\mathbb
D_d^{-1/2}Z^{(d)}\|^\gamma\geqslant d^{\gamma/2}$.

On the other hand, in the general case, the application of
\eqref{bbn} may not lead to a loss of precision, while the
statement of Theorem \ref{th9} is simpler than that of
Theorem~\ref{th6}. It involves only the moment ${\mathbf E}\,
\|Z\|^\gamma$ and the eigenvalues of the covariance operator
of~$\mathbb D$ of the vector~$Z$. An intermediate situation is
possible, where inequality \eqref{bn9} is not valid, but the
statement of Theorem \ref{th6} is still stronger than that of
Theorem~\ref{th9}.

The proofs of Theorems \ref{th6} and \ref{th9} are based on the
method of finite dimensional approximation, related to the method
applied for estimating the accuracy of approximation in the CLT in
infinite dimensional spaces (see., for instance, the survey
\cite{78}).

Applying Chebyshev's inequality, we see that, under the
assumptions of Theorem~\ref{th3}, we have
\begin{equation}
{\mathbf P} \bigl\{\Delta_n (X, Y)\geqslant x\bigr\}
\ll_{\gamma,d} \bigl( \sigma_{1}/\sigma_{d}\bigr) ^\gamma \, n\,
{\mathbf E}\, \|Z\|^\gamma/x^\gamma\label{1.13}
\end{equation}
for all $x>0$ and all $n=1, 2, \ldots$. Clearly, the statement of
Theorem~\ref{th3} is stronger than  \eqref{1.13}. A construction
for which  \eqref{1.13} is valid for $d=1$ for fixed $n$ and
$x=O\big(\sqrt n\,\log n\big)$  with constants, depending on
$\gamma$ and $\mathcal L(Z)$ only,  was proposed by Koml\'os,
Major, and Tusn\'ady (KMT)~\cite{7}, see also Borovkov~\cite{1}
and Major~\cite{9} in the case $2<\gamma\leqslant3$. Then
Sakhanenko \cite{12} proved Theorem~\ref{th2}, which ensures the
validity of the one-dimensional version of inequality~\eqref{1.13}
for all $x$ on the same probability space. Einmahl~\cite{4}
obtained a multidimensional version of the KMT result without
restrictions on the values of~$x$.

Previously, the estimates for the rate of strong approximation in
infinite dimensional spaces appeared, for example, in \cite{21},
\cite{31}, \cite{41} and~\cite{S2}. The closest to the subject of
this paper is the following infinite dimensional result of
Sakhanenko~\cite{S2}.

\begin{theorem} \label{th8} Let $Z$
be an $ \mathbf H$-valued random vector with ${\mathbf E}\,Z=0$
and ${\mathbf E}\,\left\Vert  Z\right\Vert ^\gamma <\infty$, for
some $\gamma$ with $2\leqslant\gamma\leqslant3$. Then, for any
fixed $x>0$, there exists a construction such that
\begin{equation}
{\mathbf P} \bigl\{\Delta_n^\infty (X, Y)\geqslant x\bigr\} \ll
n\, {\mathbf E}\, \|Z\|^\gamma/x^\gamma\quad\text{for all }
n,\label{1.13y}
\end{equation}
where
\begin{equation}
\Delta_n^\infty (X,Y)=\max_{1\leqslant s\leqslant n}\;\Bigl\|\,\sum_{j=1}^sX_j-
\sum_{j=1}^sY_j\,\Bigr\|_\infty.  \label{1.1w}
\end{equation}
\end{theorem}\medskip

Theorem~\ref{th8} is formulated for fixed $x$. This means that the
probability space depends on this~$x$. Furthermore, in the
statement of Theorem~\ref{th8}, the quantity  $\Delta_n (X, Y)$ is
replaced by $\Delta_n^\infty (X, Y)$, which is (in general)
essentially smaller than $\Delta_n (X, Y)$. On the other hand,
inequality~\eqref{1.13y} looks almost as inequality \eqref{1.13}
for~$2\leqslant\gamma\leqslant3$. We should note that Sakhanenko
\cite{S2} obtained substantially more general results in
comparison with Theorem~\ref{th8}. They are proved for
non-identically distributed depending summands, forming, for
example, infinite dimensional martingales.

The following theorem yields a lower bound under the assumptions
of Theorems~\ref{th6} and \ref{th9}.

\begin{theorem}\label{th11} Let positive numbers
$\sigma_m^2$, $m=1,2,\dots$, satisfy the relations
 \begin{equation}\label{in441}
\sigma_1^2\geqslant\sigma_2^2\geqslant\dots\geqslant\sigma_m^2\geqslant\dots,\quad
\hbox{�} \quad\sum_{m=1}^\infty\sigma_m^2<\infty.
\end{equation}
Let  $n$ be a fixed positive integer, and $\lambda>0$ with
$\sigma_1^2\leqslant\lambda^{2}$. Denote
\begin{equation}\label{16w}k
=\min\{m:n\,\sigma_m^2<\lambda^{2}\}-1.\end{equation} Then there
exists an $ \mathbf H$-valued random vector $Z=(Z_1,Z_2,\dots)$,
satisfying \eqref{in443}--\eqref{in45} and such that ${{\mathbf
E}\,\left\|Z\right\|^\gamma<\infty}$, for all
${\gamma\geqslant0}$, and for any construction we have the lower
bound
\begin{equation}
{\mathbf E}\, \bigl(\Delta_n (X, Y)\bigr)^\gamma\gg_\gamma
{\mathbf E}\, \bigl(\Delta_n (X^{(k)}, Y^{(k)})\bigr)^\gamma+(n
B_k ^2)^{\gamma/2}.\label{1.10mm}
\end{equation} Meanwhile, the first term in the right-hand side of \eqref{1.10mm}
is assumed to be zero if $k =0$.\end{theorem}

\begin{remark}\rm
Finding a lower bound for ${\mathbf E}\, \bigl(\Delta_n (X^{(k)},
Y^{(k)})\bigr)^\gamma$ is a separate problem. Note, however, that
the vector $Z$ from the proof of Theorem~\ref{th11} satisfies the
rough bound ${\mathbf E}\, \bigl(\Delta_n (X^{(k)},
Y^{(k)})\bigr)^\gamma\gg_{\gamma}(\lambda^{2} k)^{\gamma/2}$,
since it has a lattice distribution.
\end{remark}

 The presence of the
quantity $(nB_k^2)^{\gamma/2}$ in the right-hand side of
\eqref{1.10mm} confirms that the appearance
 of the summand $(nB_d^2)^{\gamma/2}$ in  \eqref{1.10r} and~\eqref{1.10hr} is natural.
 It becomes clear when we compare inequality \eqref{1.10mm}
with the intermediate inequality \eqref{in6}.

In  Section~3, we consider Examples 1--4, showing, in particular,
that for many distributions Theorem~\ref{th9} yields estimates,
which are stronger than the estimates of Theorem~\ref{th8}.
Moreover, in Example 5, we verify that, if the sequence of
eigenvalues $\sigma_{m}^{2}$ decreases slowly, then Theorems~
\ref{th6} and~\ref{th9} provide estimates which are optimal in
order.

\section{Proofs}

We shall need the following Lemmas~\ref{lr1}--\ref{lms}.\medskip

\begin{lemma}\label{lr1} Let $\xi _1,\dots ,\xi _n$
be
independent random vectors which have mean zero and assume values
in $\mathbf  H$. Then
\begin{equation}
{\mathbf E}\, \Bigl\| \sum\limits_{j=1}^n \xi_j\Bigr \|^\gamma
\ll_\gamma \sum\limits_{j=1}^n {\mathbf E}\, \|\xi_j\|^\gamma +
\Bigl( \sum\limits_{j=1}^n {\mathbf E}\, \|\xi_j\|^2\Bigr)^{\gamma
/2}, \quad \hbox{for  $\gamma\geqslant 2$} . \label{1.15}
\end{equation}
\end{lemma}

This multidimensional version of the Rosenthal inequality follows
easily from a result of de Acosta~\cite{2}. In the i.i.d. case,
the second summand in the right-hand side of~\eqref{1.15} grows
faster than the first term as $n\to\infty$. Theorems~\ref{th2}
and~\ref{th3} show that this growth corresponds to the growth of
moments of sums of Gaussian approximating vectors.

The next lemma is proved by Rosenthal~\cite{10}, see also Johnson,
Schechtman and Zinn~\cite{6}.

\begin{lemma}\label{lr2} Let \,$\xi _{1},\dots ,\xi _{n}$ \,be
independent random variables which are non-negative with
probability one. Then
$$
{\mathbf E}\, \Bigl( \sum\limits_{j=1}^n \xi_j\Bigr )^\gamma
\ll_\gamma \sum\limits_{j=1}^n {\mathbf E}\, \xi_j^\gamma + \Bigl(
\sum\limits_{j=1}^n {\mathbf E}\, \xi_j\Bigr)^\gamma  \quad
\hbox{for }\gamma\geqslant 1.
$$
\end{lemma}

The following Lemma~\ref{lms} is proved by Montgomery-Smith \cite{MS}. It is a particular case of Theorem 1.1.5 from
the monograph of de la Pe$\tilde{\text{n}}$a and Gin\'e~\cite{3}.

\begin{lemma} \label{lms} Let \,$\xi _{1},\dots ,\xi _{n}$ \,be i.i.d.
random vectors with values in $\mathbf H$. Then
$$
{\mathbf P} \Bigl\{\max\limits_{1\leqslant s\leqslant n}\Bigl\|
\sum\limits_{j=1}^s \xi_j\Bigr \|>x\Bigr\} \leqslant 9\,{\mathbf P}
\Bigl\{\Bigl\| \sum\limits_{j=1}^n \xi_j\Bigr \|>x/30\Bigr\} \quad\hbox{for all }x\geqslant0.
$$
\end{lemma}

Coupled with the well-known equality
$$
{\mathbf E}\, \left|\eta\right|^\gamma=\gamma\int_0^\infty
x^{\gamma-1}\, \mathbf P\bigl\{\left|\eta\right|> x\bigr\}\,dx,
\quad\gamma>0,
$$ which is
valid for any random variable~$\eta$, Lemma~\ref{lms} allows us to
estimate the moments
\begin{equation}
{\mathbf E}\, \max\limits_{1\leqslant s\leqslant n}\Bigl\| \sum\limits_{j=1}^s
\xi_j\Bigr \|^\gamma\ll_\gamma{\mathbf E}\, \Bigl\| \sum\limits_{j=1}^n
\xi_j\Bigr \|^\gamma,\quad \gamma>0,\label{1.16}
\end{equation}
in the case of i.i.d. random vectors $\xi_1,\dots ,\xi _n$.
\medskip

\noindent {\bf Proof of Theorem \ref{th6}.}
It is not difficult to understand that for any construction we have
\begin{eqnarray}
{\mathbf E}\, \bigl(\Delta_n (X, Y)\bigr)^\gamma\ll_\gamma &&
\hskip-.5cm{\mathbf E}\, \bigl(\Delta_n (X^{(d)}, Y^{(d)})\bigr)^\gamma
\nonumber\\ &&\hskip -1cm +\,{\mathbf E}\,
\max_{1\leqslant s\leqslant n} \Bigl\| \sum\limits_{j=1}^s X_j^{[d]}\Bigr \|^\gamma
+{\mathbf E}\,\max_{1\leqslant s\leqslant n} \Bigl\| \sum\limits_{j=1}^s Y_j^{[d]}\Bigr \|^\gamma .\label{in9}
\end{eqnarray}
Using \eqref{1.16} and \eqref{in9}, we obtain
\begin{equation}\label{in3}
{\mathbf E}\, \bigl(\Delta_n (X, Y)\bigr)^\gamma\ll_\gamma {\mathbf E}\,
\bigl(\Delta_n (X^{(d)}, Y^{(d)})\bigr)^\gamma+
{\mathbf E}\, \Bigl\| \sum\limits_{j=1}^n X_j^{[d]}\Bigr \|^\gamma +
{\mathbf E}\, \Bigl\| \sum\limits_{j=1}^n Y_j^{[d]}\Bigr \|^\gamma .
\end{equation}

Lemma \ref{lr1} together with $\mathcal L(X_j)=\mathcal L(Z)$ and
$\hbox{\rm cov}\,Y_j=\hbox{\rm cov}\,X_j=\hbox{\rm cov}\,Z$ yields
\begin{equation}\label{in4}
{\mathbf E}\, \Bigl\| \sum\limits_{j=1}^n X_j^{[d]}\Bigr \|^\gamma \ll_\gamma
n\, {\mathbf E}\, \bigl\|Z^{[d]}\bigr\|^\gamma + \bigl(n\,
 {\mathbf E}\, \bigl\|Z^{[d]}\bigr\|^2\bigr)^{\gamma/2} ,
\end{equation}
and
\begin{equation}\label{in5}
{\mathbf E}\, \Bigl\| \sum\limits_{j=1}^n Y_j^{[d]}\Bigr \|^\gamma
\ll_\gamma \bigl( n\, {\mathbf E}\,
\bigl\|Z^{[d]}\bigr\|^2\bigr)^{\gamma/2}.
\end{equation}
Inequalities \eqref{in3}--\eqref{in5}  imply that
\begin{equation}\label{in6}
{\mathbf E}\, \bigl(\Delta_n (X, Y)\bigr)^\gamma\ll_\gamma
{\mathbf E}\, \bigl(\Delta_n (X^{(d)}, Y^{(d)})\bigr)^\gamma+
n \,{\mathbf E}\, \bigl\|Z^{[d]}\bigr\|^\gamma +  (n\,B_d^2)^{\gamma/2} .
\end{equation}

It is easy to show that condition \eqref{1.10ru} implies that the
$d$-dimensional vector~$Z^{(d)}$ satisfies condition
\eqref{1.10bt} of Theorem~\ref{th4}. Applying that theorem, we see
that from \eqref{1.10b} and from the well-known  Berkes--Philipp
Lemma~\cite{BF} it follows, that there exists a construction such
that
\begin{equation}
{\mathbf E}\, \bigl(\Delta_n (X^{(d)},
Y^{(d)})\bigr)^\gamma\ll_{\gamma, \psi}d^{\psi\gamma}\, n\,
\sigma_{1}^\gamma \, {\mathbf E}\, \bigl\|\mathbb
D_d^{-1/2}Z^{(d)}\bigr\|^\gamma.\label{1.10aa}
\end{equation}
Using \eqref{in6} and \eqref{1.10aa}, we obtain the  statement of
Theorem~\ref{th6}.
\medskip

\noindent {\bf Proof of Theorem \ref{th11}.} Let $Z_j$
(coordinates of the vector $Z$) be independent random variables,
taking values $-\lambda$,  0 and $\lambda$ with probabilities
\begin{equation}
{\bf P}\{ Z_m=\pm\lambda\}=\sigma_m^2/2\lambda^{2},\quad{\bf P}\{
Z_m=0\}=1-\sigma_m^2/\lambda^{2},\quad m=1,2, \dots.\label{1.11p}
\end{equation}
 With the help of Lemma \ref{lr2} it
is not difficult to show that ${\mathbf
E}\,\left\|Z\right\|^\gamma<\infty$, for all $\gamma\geqslant0$.
Assume that we have constructed a sequence of independent random
vectors $X_1,\dots ,X_n$  and a corresponding sequence of
independent Gaussian random vectors  $Y_1,\dots ,Y_n$ such that $
\mathcal L(X_j)=\mathcal L(Z)$, ${\mathbf E}\,Y_j=0$, $\hbox{\rm
cov}\,Y_j=\hbox{\rm cov}\,X_j$, $j=1,\dots ,n$.

Then the coordinates of the vectors \,$X_j$ \, (namely \{$X_{jm}$,
$j=1,2, \dots,n$, $m=1,2, \dots$\}) are jointly independent random
variables with distributions ${\mathcal L}(Z_m)$, while the
coordinates of the vectors \,$Y_j$ \, (namely \{$Y_{jm}$, $j=1,2,
\dots,n$,  $m=1,2, \dots$\}) are jointly independent Gaussian
random variables with mean zero and variances~$\sigma_m^2$. Set
\begin{equation}
S_{nm}=\sum_{j=1}^nX_{jm} ,\quad T_{nm}=\sum_{j=1}^nY_{jm}, \quad
m=1,2, \dots .\label{1.1k}
\end{equation}
It is clear that ${\hbox{\rm Var}}\, S_{nm}={\hbox{\rm Var}}\,
T_{nm}=n\,\sigma_m^2$, for $m=1, 2, \dots$, and
\begin{equation}
\Delta_n (X, Y)\geqslant\max\big\{\Delta_n (X^{(k)}, Y^{(k)}),\
\Delta_n (X^{[k]}, Y^{[k]})\big\}.\label{1.19h}
\end{equation}
Obviously,
\begin{equation}
\Delta_n (X^{[k]},
Y^{[k]})\geqslant\Bigl\|\,\sum_{j=1}^nX_j^{[k]}-
\sum_{j=1}^nY_j^{[k]}\,\Bigr\|,\label{1.10h}
\end{equation}
while \begin{equation}\label{16k} \Bigl\|\,\sum_{j=1}^nX_j^{[k]}-
\sum_{j=1}^nY_j^{[k]}\,\Bigr\|^2 =
\sum_{m=k+1}^\infty\left|S_{nm}- T_{nm}\right|^2  .
\end{equation}

If $m>k $, then \begin{equation}\label{1k} \left|S_{nm}-
T_{nm}\right|\geqslant\eta_{nm},\quad\hbox{where
}\eta_{nm}=\left|T_{nm}\right|\,{\mathbf1}
\{\left|T_{nm}\right|\leqslant \lambda/2\},
\end{equation}
since the random variables $ S_{nm}$ take only values which are
multiples of~$\lambda$.
 Put
\begin{equation}
U_{nk }=  \sum_{m=k +1}^\infty\eta_{nm}^2 .\label{11k}
\end{equation}

For fixed $n$, the set $\{\eta_{nm}\}$ is a collection of jointly
independent random variables. According to \eqref{16w},
\eqref{1.1k} and \eqref{1k}, for $m>k $,
\begin{equation}\label{15s}
{\mathbf E}\,(\eta_{nm}^2)\asymp n\,\sigma_m^2\quad\hbox{and}
\quad\hbox{Var}(\eta_{nm}^2)\asymp n^2\sigma_m^4.
\end{equation}
Denote $a={\mathbf E}\,  U_{nk}$ and $b={\hbox{\rm Var}}\,
U_{nk}$. Note that by relations \eqref{16w}, \eqref{1k},
\eqref{11k} and \eqref{15s},
\begin{equation}
\label{17gw} a=  \sum_{m=k +1}^\infty{\mathbf
E}\,(\eta_{nm}^2)\asymp n\,B_k ^2 \quad\text{and}\quad b=
\sum_{m=k +1}^\infty\hbox{Var}(\eta_{nm}^2)\ll a^2 ,
\end{equation}
where the quantity $B_k ^2$ is defined by formula \eqref{in44}.
According to inequality (7.5) from Feller \cite{Fe}, p. 180,
\begin{equation}\label{17g}
{\mathbf P}\big\{  U_{nk }-a<-t
\big\}\leqslant\frac{b}{b+t^2}=1-\frac{t^2}{b+t^2},\quad\hbox{for
all }t\geqslant0.
\end{equation}
Applying \eqref{17g} for $t=a/2$ and relations \eqref{17gw}, it is
easy to show that
\begin{equation}\label{18j}
{\mathbf P}\big\{  U_{nk }\geqslant a/2\big\}\gg1 .
\end{equation}
Therefore, relations \eqref{1k}, \eqref{11k}
 and \eqref{18j} yield
\begin{equation}\label{18kk}
{\mathbf P}\Bigl\{\sum_{m=k +1}^\infty\left|S_{nm}-
T_{nm}\right|^2 \geqslant a/2\Bigr\}\gg1.
\end{equation}
>From \eqref{1.10h}, \eqref{16k} and \eqref{18kk} we obtain
\begin{equation}\label{18p}
{\mathbf P}\bigl\{\bigl(\Delta_n (X^{[k]}, Y^{[k]})\bigr)^2
\geqslant a/2\bigr\}\gg1
\end{equation}
and, hence,
\begin{equation}
{\mathbf E}\, \bigl(\Delta_n (X^{[k]},
Y^{[k]})\bigr)^\gamma\gg_\gamma a^{\gamma/2}\asymp_\gamma(n B_k
^2)^{\gamma/2}.\label{1.10m}
\end{equation}
Finally \eqref{1.19h} and \eqref{1.10m} imply the lower bound
 \eqref{1.10mm}.
\medskip

\section{Examples}

 In Examples 1--5 we compare the estimates which follows from Theorem \ref{th6}
 when condition \eqref{bn9} is satisfied to bounds of Theorem~\ref{th9}
 for concrete sequences of eigenvalues of the covariance
operator of the vector~$Z$.
\medskip

{\bf Example 1.} Let $\sigma_m^2=\exp\{-\alpha\,m^\beta\}$,
$m=1,2,\dots$, where $\alpha,\beta>0$. Assume that $n$ is so large
that
\begin{equation}\label{66w}d=\max\{m:n^{2/\gamma}\,(\log ^\ast
n)^{2\psi/\beta}/\sigma_{m}^2 \, < n   \sigma_{m}^2
\}\geqslant1.\end{equation} Then it is clear that
\begin{equation}\label{67w}d\asymp_{\alpha,\beta}(\log^*n )^{1/\beta}
\end{equation}
and
\begin{equation}\label{68w}\sigma_{d+1}^4 \leqslant n^{-1+2/\gamma}\,(\log ^\ast
n)^{2\psi/\beta}\leqslant \sigma_{d}^4.\end{equation} Thus, for
sufficiently large $n$, the right-hand side of inequality
\eqref{1.10hr} admits the upper bound
\begin{equation}\label{69h} d^{\psi\gamma}\bigl(
\sigma_{1}/\sigma_{d}\bigr)^\gamma \, n\, {\mathbf E}\,
\|Z\|^\gamma +
 n ^{\gamma/2}  B_d^\gamma\ll_{\alpha,\beta, \gamma} n^{(2+\gamma)/4}\,(\log ^\ast
n)^{\psi\gamma/2\beta}\, {\mathbf E}\, \|Z\|^\gamma.\end{equation}
Using relations \eqref{67w} and \eqref{68w}, it is not difficult
to verify that, for sufficiently large $n$, condition
\eqref{1.10hu} is satisfied and, hence, the statement of Theorem
\ref{th9} is valid with the estimate
\begin{equation}\label{59h} {\mathbf E}\, \bigl(\Delta_n (X, Y)\bigr)^\gamma
\ll_{\alpha,\beta, \gamma,\psi} n^{(2+\gamma)/4}\,(\log ^\ast
n)^{\psi\gamma/2\beta}\, {\mathbf E}\, \|Z\|^\gamma.\end{equation}
The right-hand side of inequality \eqref{59h} grows slower than
$n^{\gamma/2}$ (the order of the trivial estimate which follows
from Lemma~\ref{lr1} and inequality~\eqref{1.16}). Therefore,
inequality \eqref{59h} is a meaningful estimate of the rate of
approximation in the infinite dimensional invariance principle. In
particular, using Lyapunov's inequality ${\mathbf E}\,
\Delta^3\leqslant ({\mathbf E}\, \Delta^\gamma)^{3/\gamma}$, we
obtain that, for $\gamma>3$,
\begin{equation}\label{49h} {\mathbf E}\, \bigl(\Delta_n (X, Y)\bigr)^3
\ll_{\alpha,\beta, \gamma,\psi} n^{3(2+\gamma)/4\gamma}\,
({\mathbf E}\, \|Z\|^\gamma)^{3/\gamma}.\end{equation} For
$\gamma>6$, the order of inequality \eqref{49h} with respect to
$n$ is better than the order of estimate~\eqref{1.13y}.
\medskip

{\bf Example 2.} Suppose now, under the assumptions of Example 1,
that \eqref{bn9} holds. Assume that $n$ is so large that
\begin{equation}\label{g66w}d=\min\{m: n   B_{m}^2<1
\}\geqslant1.\end{equation} It is clear that then relation
\eqref{67w} is still satisfied. Thus, for sufficiently large $n$,
the right-hand side of inequality \eqref{1.10r} admits the upper
bound
\begin{eqnarray}\label{g69h}  d^{\psi\gamma}\, n\, \sigma_{1}^\gamma \, {\mathbf E}\,
\bigl\|\mathbb
D_d^{-1/2}Z^{(d)}\bigr\|^\gamma\hskip-1.5cm&+\hskip-1.5cm& \;n
\,{\mathbf E}\, \bigl\|Z^{[d]}\bigr\|^\gamma+
 (n\,B_d^2)^{\gamma/2}\nonumber\\ &\ll_{\alpha,\beta, \gamma, K}& n\,(\log ^\ast
n)^{(2\psi+1)\gamma/2\beta}\, {\mathbf E}\,
\|Z\|^\gamma.\end{eqnarray} Using relation \eqref{67w}, it is also
not difficult to verify that, for sufficiently large $n$,
condition \eqref{1.10ru} is satisfied and, hence, the statement of
Theorem~\ref{th6} is valid with the estimate
\begin{equation}\label{g59h} {\mathbf E}\, \bigl(\Delta_n (X, Y)\bigr)^\gamma
\ll_{\alpha,\beta, \gamma,\psi, K} n\,(\log ^\ast
n)^{(2\psi+1)\gamma/2\beta}\, {\mathbf E}\,
\|Z\|^\gamma,\end{equation} which is considerably stronger than
\eqref{59h} and is close to the finite-dimensioval estimate
\eqref{1.10} of Theorem \ref{th3}.
\medskip

{\bf Example 3.} Let $\sigma_m^2=m^{-b}$, $m=1,2,\dots$, where
$b>1$. Choose
\begin{equation}\label{f66w}d=\max\{m:n^{2/\gamma}\,m^{2\psi}/\sigma_{m}^2 \,
< n  m\sigma_{m}^2 \}.\end{equation} It is clear that then
$d\geqslant1$ and
\begin{equation}\label{f67w}d^{b-1}\asymp_b n^{r(\gamma-2)/\gamma}
,\quad\hbox{where } r=\frac{b-1}{2b-1+2\psi}.
\end{equation}
Therefore, if \eqref{1.10hu} is fulfilled, then the right-hand
side of inequality \eqref{1.10hr} admits the upper bound
\begin{equation}\label{f69h} d^{\psi\gamma}\bigl(
\sigma_{1}/\sigma_{d}\bigr)^\gamma \, n\, {\mathbf E}\,
\|Z\|^\gamma +
 n ^{\gamma/2}  B_d^\gamma\ll_{b, \gamma}
 n^{(\gamma-r(\gamma-2))/2}\, {\mathbf E}\, \|Z\|^\gamma.\end{equation}
Using \eqref{f67w}, it is not difficult to verify that, for
sufficiently large $n$, condition \eqref{1.10hu} is satisfied
provided that $ \gamma<2\,(b-1+2\psi)$. In this case the statement
of Theorem \ref{th9} is valid with the estimate
\begin{equation}\label{f59h} {\mathbf E}\, \bigl(\Delta_n (X, Y)\bigr)^\gamma
\ll_{\psi,b, \gamma} n^{(\gamma-r(\gamma-2))/2}\,{\mathbf E}\,
\|Z\|^\gamma.\end{equation} Using Lyapunov's inequality by analogy
with Example 1, we obtain that for $\gamma>3$
\begin{equation}\label{f49h} {\mathbf E}\, \bigl(\Delta_n (X, Y)\bigr)^3
\ll_{\psi,b, \gamma} n^{3(\gamma-r(\gamma-2))/2\gamma}\, ({\mathbf
E}\, \|Z\|^\gamma)^{3/\gamma}.\end{equation} For $3r-1>0$ and
$\gamma>6r/(3r-1)$, the order of inequality \eqref{f49h} with
respect to $n$ is better than the order of estimate~\eqref{1.13y}.

If condition \eqref{1.10hu} is not fulfilled for $d$ defined by
 \eqref{f66w}, one should decrease
$d$ choosing
\begin{equation}\label{f49hh}
d=\max\bigl\{m:C(\gamma)\, m^{\gamma/2}(\log^* m)^{\gamma+1}\,
\bigl(\mathbf E\left\|Z\right\|^\gamma\bigr)^{2/\gamma} \leqslant
n^{1-2/\gamma}\,\sigma_m^2 \bigr\}.
\end{equation}
 It is clear that then, for
sufficiently large $n$, we have
\begin{equation}\label{f6y}d^{b-1}\asymp_{b, \gamma,\lambda} n^{\delta(\gamma-2)/\gamma}
(\log^*n )^{-\delta(\gamma+1)},\quad\hbox{��� }
\delta=\frac{2(b-1)}{2b+\gamma}, \ \lambda={\mathbf E}\,
\|Z\|^\gamma.
\end{equation}
In this case the statement of Theorem \ref{th9} is valid with the
estimate
\begin{equation}\label{f55h} {\mathbf E}\, \bigl(\Delta_n (X, Y)\bigr)^\gamma
\ll_{\psi,b, \gamma, \lambda}
n^{(\gamma-\delta(\gamma-2))/2}\,(\log ^\ast
n)^{\delta\gamma(\gamma+1)/2},\end{equation} which must be weaker
in order in comparison with \eqref{f59h}. Thus, in the general
case, for sufficiently large $n$, there exists a construction such
that
\begin{eqnarray}\label{f77} &&\hskip-1cm{\mathbf E}\,
\bigl(\Delta_n (X, Y)\bigr)^\gamma \nonumber \\
 &&\ll_{\psi,b, \gamma, \lambda}\max\big\{
n^{(\gamma-\delta(\gamma-2))/2}\,(\log ^\ast
n)^{\delta\gamma(\gamma+1)/2},\,n^{(\gamma-r(\gamma-2))/2}\big\}.\end{eqnarray}
\medskip

{\bf Example 4.} Suppose now, under the assumptions of Example 3,
that \eqref{bn9} holds. Choose
\begin{equation}\label{gf66w}d=\min\{m:n^{2/\gamma}\,m^{2\psi+1}
< n  m\sigma_{m}^2 \}.\end{equation} It is clear that then
$d\geqslant1$ and
\begin{equation}\label{gf67w}d^{b-1}\asymp_b n^{\rho (\gamma-2)/\gamma}
,\quad\hbox{where } \rho =\frac{b-1}{b+2\psi}.
\end{equation}
Therefore, if \eqref{1.10ru} is fulfilled, then the right-hand
side of inequality \eqref{1.10r} admits the upper bound
\begin{eqnarray}\label{gf69h}&& \hskip-2cm d^{\psi\gamma}\, n\, \sigma_{1}^\gamma \, {\mathbf E}\,
\bigl\|\mathbb D_d^{-1/2}Z^{(d)}\bigr\|^\gamma+ \;n \,{\mathbf
E}\, \bigl\|Z^{[d]}\bigr\|^\gamma+
 (n\,B_d^2)^{\gamma/2}\nonumber \\
 &&\hskip3cm\ll_{b, \gamma, K}
 n^{(\gamma-\rho (\gamma-2))/2}\, {\mathbf E}\, \|Z\|^\gamma.
 \end{eqnarray}
Using \eqref{gf67w}, it is also not difficult to verify that, for
sufficiently large $n$, condition \eqref{1.10ru} is satisfied
provided that $ \gamma<2\,(b-1+2\psi)$. In this case the statement
of Theorem \ref{th6} is valid with the estimate
\begin{equation}\label{gf59h} {\mathbf E}\, \bigl(\Delta_n (X, Y)\bigr)^\gamma
\ll_{\psi,b, \gamma, K} n^{(\gamma-\rho (\gamma-2))/2}\,{\mathbf
E}\, \|Z\|^\gamma.\end{equation}

If condition \eqref{1.10ru} is not fulfilled, then
 one should choose $d$ for applying
Theorem \ref{th6} not by formula \eqref{gf66w}, but by relation
\begin{equation}\label{gf49hh}
d=\max\bigl\{m: C(\gamma)\,K^{2/\gamma}\, m^{1+\gamma/2}(\log^*
m)^{\gamma+1} \leqslant n^{1-2/\gamma}\bigr\}.
\end{equation}
 It is clear that, for sufficiently large $n$,
\begin{equation}\label{gf6y}d^{b-1}\asymp_{b, \gamma,K} n^{\mu(\gamma-2)/\gamma}
(\log^*n )^{-\mu(\gamma+1)},\quad\hbox{where }
\mu=\frac{2(b-1)}{2+\gamma}.
\end{equation}
In this case the statement of Theorem \ref{th6} is valid with the
estimate
\begin{equation}\label{gf55h} {\mathbf E}\, \bigl(\Delta_n (X, Y)\bigr)^\gamma
\ll_{\psi,b, \gamma, K} n^{(\gamma-\mu(\gamma-2))/2}\,(\log ^\ast
n)^{\mu\gamma(\gamma+1)/2},\end{equation} which is weaker in order
in comparison with \eqref{gf59h}. In the general case, for
sufficiently large $n$, there exists a construction such that
\begin{eqnarray}\label{f77} {\mathbf E}\,
\bigl(\Delta_n (X, Y)\bigr)^\gamma \ll_{\psi,b, \gamma,
K}{\hskip2cm}&&\nonumber\\ \max\big\{
n^{(\gamma-\mu(\gamma-2))/2}\,(\log ^\ast
n)^{\mu\gamma(\gamma+1)/2},&&\,\hskip-0.5cm
n^{(\gamma-\rho(\gamma-2))/2}\big\}\,{\mathbf E}\,
\|Z\|^\gamma.\end{eqnarray}
\medskip

{\bf Example 5.} Let $\sigma_m^2=1/m\,(\log^*m)^{1+\tau}$, for
$m=1,2,\dots$, where $\tau>0$. Denote by $ \lceil x \rceil$
 the integer part of a number~$x$. Choose
\begin{equation}\label{yg66w}d=\lceil n^{\varepsilon}\rceil,
\quad\hbox{where }
\varepsilon=\frac{\gamma-2}{\gamma\,(\gamma+22)}.
\end{equation} It is clear that then
\begin{equation}\label{yg67w}B_d^{2}\asymp_{\tau}\frac 1{(\log^*d)^{\tau}}
\asymp_{\gamma,\tau}\frac 1{(\log^*n)^{\tau}}.
\end{equation}
 Using relations \eqref{yg66w} and
\eqref{yg67w}, it is not difficult to verify that, for
sufficiently large $n$, condition \eqref{1.10hu} is satisfied and
the statement of Theorem~\ref{th9} is valid with the estimate
\begin{equation}\label{yg59h} {\mathbf E}\, \bigl(\Delta_n (X, Y)\bigr)^\gamma
\ll_{\gamma,\tau} (n/(\log ^\ast
n)^{\tau})^{\gamma/2}.\end{equation}
\medskip

Let us compare the upper bounds obtained in Examples 1, 3 and 5
using Theorem~\ref{th9} and the lower bound
\begin{equation}
{\mathbf E}\, \bigl(\Delta_n (X, Y)\bigr)^\gamma\gg_\gamma (n B_k
^2)^{\gamma/2},\label{3mm}
\end{equation}
which follows from \eqref{1.10mm}.

In Example 1 the lower bound \eqref{3mm} is far from the upper
bound \eqref{59h}. It is not difficult to calculate that, in
Example~3, the positive integer $k$, defined by \eqref{16w},
satisfies the relations $k\asymp_{b,\lambda} n^{1/b}$,
$B_k^2\asymp_{b,\lambda} n^{(1-b)/b}$, while the lower bound
\eqref{3mm} is of order $O\bigl(n^{\gamma/2b}\bigr)$. This shows,
that the order of upper bounds should be expected to grow when
$\gamma$ increases. Notice that, for large values of $\gamma$ and
$b$, the order of upper bounds is close to $n^{\gamma/4}$.

 For
relatively small values of $\gamma$ and $b$, the orders of
estimates depend essentially on~$\psi$, which is involved in the
bounds due to sufficiently large powers of dimension $d$ in the
estimates of Theorems~\ref{th3} and~\ref{th4}. Possible
improvements of Theorems \ref{th3} and~\ref{th4} have to improve
the order of the upper bounds in Examples~1--4 and in
Theorems~\ref{th6} and ~\ref{th9}.

In Example 5, it is easy to verify that $k\asymp_{\tau} n/(\log^*
n)^{1+\tau}$. Thus, the upper and lower bounds are of the same
order $O\bigl((n/(\log ^\ast n)^{\tau})^{\gamma/2}\bigr)$, and
Theorem \ref{th9} provides the correct order for the rate of
approximation. The same is true if the variances of coordinates
$\sigma_m^2$ are decreasing slower than in Example~5. Therefore,
the order of estimates could be made arbitrarily close to the
trivial order $O\bigl(n^{\gamma/2}\bigr)$.

The authors are grateful to a referee for a series of useful
remarks which enable us to improve exposition substantially.

 \end{document}